# *A Systematic Literature Review*: Ethnomathematics in Geometry

**Ratu Sarah Fauziah Iskandar, Natanael Karjanto, Yaya Sukjaya Kusumah, Iden Rainal Ihsan**

[a]Department of Mathematics Education, Universitas Muhammadiyah Tangerang,
Jalan Perintis Kemerdekaan I/33, Tangerang, Indonesia
[b]Department of Mathematics, Sungkyunkwan University, Natural Science Campus 2066 Seobu-ro,
Jangan-gu, Suwon 16419, Gyeonggi-do, Republic of Korea
[c]Department of Mathematics Education, Universitas Pendidikan Indonesia,
Jl. Dr. Setiabudhi No. 229, Bandung, Indonesia
[d]Department of Mathematics Education, Universitas Samudra,
Jln. Prof. Dr. Syarief Thayeb, Meurandeh, Langsa Lama, Langsa, Aceh, Indonesia
*Corresponding author: natanael@skku.edu

**Abstract**

This study aims to find out geometric concepts related to ethnomathematics in Indonesia. The design used in this study is to summarize, review, and analyze 24 articles with the scope of ethnomathematics discussion with journals or proceedings indexed by Scopus in the last 5 years. The year of publication of the articles is from 2016 to 2020. With the assistance of the Publish or Perish tool, the researcher was able to compile a collection of journal articles and conferences that were indexed by Scopus to finish this study. The terms to keep in mind are geometry and ethnomathematics. The researcher use PRISMA, which incorporates resources such as Scopus, to conduct the systematic review. This includes determining eligibility and exclusion criteria, as well as the steps of the review process (identification, screening, and eligibility), as well as the abstraction and analysis of data. From the results of the study obtained several cultures that use the concept of geometry including traditional houses, regional dances, regional specialties, batik motifs and fabrics / weaving, temples and mosques.

## INREODUCTION

The material taught in math lessons is relevant in daily life. School mathematics is primarily a transfer of knowledge, with students being encouraged to think critically, creatively, and reflectively about the knowledge they learn and acquire (Prahmana & Istinndaru, 2021; Risdiyanti & Prahmana, 2020)When it comes to learning mathematics, children still have trouble grasping abstract knowledge; therefore, they want a strategy that translates the notion into something more concrete in the hope that it may be utilized in everyday life (Brandt & Chernoff, 2015). Because of the frequent connections that may be found between culture and mathematics, the mathematics that is taught in schools is inextricably linked to the culture of the students (Seah & Bishop, 2003).

Mathematics and culture are two aspects of existence that just cannot be divorced from one another. Culture emerges within society as a result of the typical means through which humans adapt to their surroundings (Supiyati, Hanum & Jailani, 2019). Because every society is required to have its own culture, culture is an aspect that is very closely associated with the community. According to Gilmer (1990) teaching math without cultural context on the pretext that it is perfect, abstract, and universal is the reason for students' declining achievements and their failure.

The geometry of graphic arts from a variety of cultures is a striking example of the close connection that exists between culture and mathematics (Massarwe et al., 2010). Students often found geometry to be one of the more challenging topics. In Indonesian primary schools, the subject of geometry was greeted with a lot of challenges due to a lack of resources. For instance, the method that was employed to instruct issues related to geometry was very theoretical, and numerous abstract ideas and formulas were presented without respect to many characteristics including logic, reasoning, and comprehension. (Fauzan, Slettenhaar and Plomp, 2002).

The relationship between mathematics and everyday life has been claimed by some education experts who said that one's mathematics is influenced by their background culture because what they do is based on what they see and feel natural. Between the theory and its references to culture and the students' daily lives, the formal mathematics at the conclusion creates a significant divide. Therefore, some work is required to include cultural components into mathematics in schools. Mathematics that is related strongly to the culture mathematics model is called as ethnomathematics. In Indonesia, it has been known as ethnomathematics (Wijayanti, N.P.A, et al, 2018). The term "ethnomathematics" refers to the application, utilization, or combination of mathematical principles in the context of various cultural practices in society (D'Ambrosio U, 1985, 2016; Vasquez, 2017; Mosimege & Ismael, 2014).

The most common definition claimed that culture is an important tenet of ethnomathematics and that mathematics is not culture-free but rather culture-bound (D'Ambrosio, 2001; Rosa & Orey, 2010). This definition also suggested that mathematics is not culture-free but rather culture-bound. The ethnomathematics approach to mathematics education is seen by pupils as a solution that is both appropriate and more relevant (Rosa & Orey, 2011). In line with Supriadi et al (2019) the concept of ethnomathematics is an innovation that can be used in learning mathematics. According to Kusuma et al. (2017) culture and mathematics are closely related and in some ethnic groups, a distinctive culture is ingrained in society and may differ from other groups (Ngiza, Susanto & Lestari, 2015; Nurjanah, Mardia, I & Turmudi, 2021).

In the past few decades, numerous studies have been carried out, and all of them have reported on the challenges that students face when trying to study geometry. The fundamental cause of these issues is a disparity between the degree of instruction and the pupils' capacities for comprehending and learning new material. According to Fouze & Amit (2021) one of the leading and most successful approaches in the field of mathematics education is the ethnomathematical approach, in which instruction is based on the integration of cultural-educational elements that express mathematical values from the students' daily life.

Research conducted by Wijayanti, N.P.A, et al (2018) ethnomathematics is suggested to be the alternative for teaching and learning to elicit students to increase their geometry logic and reasoning. Based on Hendriyanto, A., Kusmayadi, T. A and L Fitriana (2020) students at the visualization level is able to visualize batik motif into geometric patterns, at the level of analysis students are able to identify what elements are known and what are needed to solve the problems, at the level of informal deductions students are able to make statements related to the relation between existing elements and able to prove the truth statement made. As a result, ethnomathematics can be used as a learning innovation to increase students' love of mathematics, motivation to learn, and creativity through their culture. From the description above, the focus of the problem in this study is how the results of the Systematic Literature Review (SLR) on the field of geometry are related to local cultural wisdom in Indonesia.

## METHODS

A Systematic Literature Review (also known as an SLR) is the method of research that was carried out. A systematic literature review (SLR) is a type of literature review that identifies, evaluates, and interprets all of the data on a research issue to respond to preset research questions. With the assistance of the Publish or Perish tool, the researcher was able to compile a collection of journal articles and conferences that were indexed by Scopus to finish this study. The terms to keep in mind are geometry and ethnomathematics. The researcher use PRISMA, which incorporates resources such as Scopus, to conduct the systematic review. This includes determining eligibility and exclusion criteria, as well as the steps of the review process (identification, screening, and eligibility), as well as the abstraction and analysis of data.

1. PRISMA
    When reviewing the articles, the PRISMA Statement, which stands for Preferred Reporting Items for Systematic Reviews and Meta-Analyses, is utilized as a set of principles. According to Sierra- Correa and Cantera Kintz (2015), it suggests three advantages, which are as follows: 1) defining obvious research questions that allow for systematic research; 2) classifying into inclusion and exclusion criteria; and 3) attempting to observe the extensive database of scientific literature in a particular period. The PRISMA Statement makes it possible to do an exhaustive search for concepts associated with the Ethnomathematics in Geometry.
2. Resources
    This review used Scopus, which is a consolidated database for journal articles. Scopus is one of the most extensive databases of abstracts and citations of literature that has been peer-reviewed. It contains over 22,800 journals from 5,000 publishers all around the world. Scopus covers a wide range of topics, including education, the social sciences, and agriculture, amongst others.
3. Eligibility and exclusion criteria
    The eligibility requirements and the restrictions are decided upon. To begin, the selection process consisted solely of published publications on ethnomathematics in geometry. Second, only scholarly journal papers on Indonesian ethnomathematics in the field of geometry. Third, a period of five years is decided upon as the term (between 2017 and 2022). This time frame is adequate for observing the progression of research and publications that are connected to it.

## RESULTS AND DISCUSSION

The systematic review process involved four stages with the assistance of the Publish or Perish tool, the researcher was able to compile a collection of journal articles and conferences that were indexed by Scopus to finish this study. The keywords are geometry and ethnomathematics. The first phase was to identify the keywords used for the searching process. At this stage, 130 duplicated articles were removed after screening. The screening process was the next step. 60 items are eligible to be examined at this stage, whereas 10 articles have been eliminated. Eligibility was the third and last stage. It was possible to read the entire article. Following careful reading and examination, a total of 36 items were deemed unsuitable for inclusion. Only 24 of the original articles were left once everything was said and done, and they were the ones that were used.

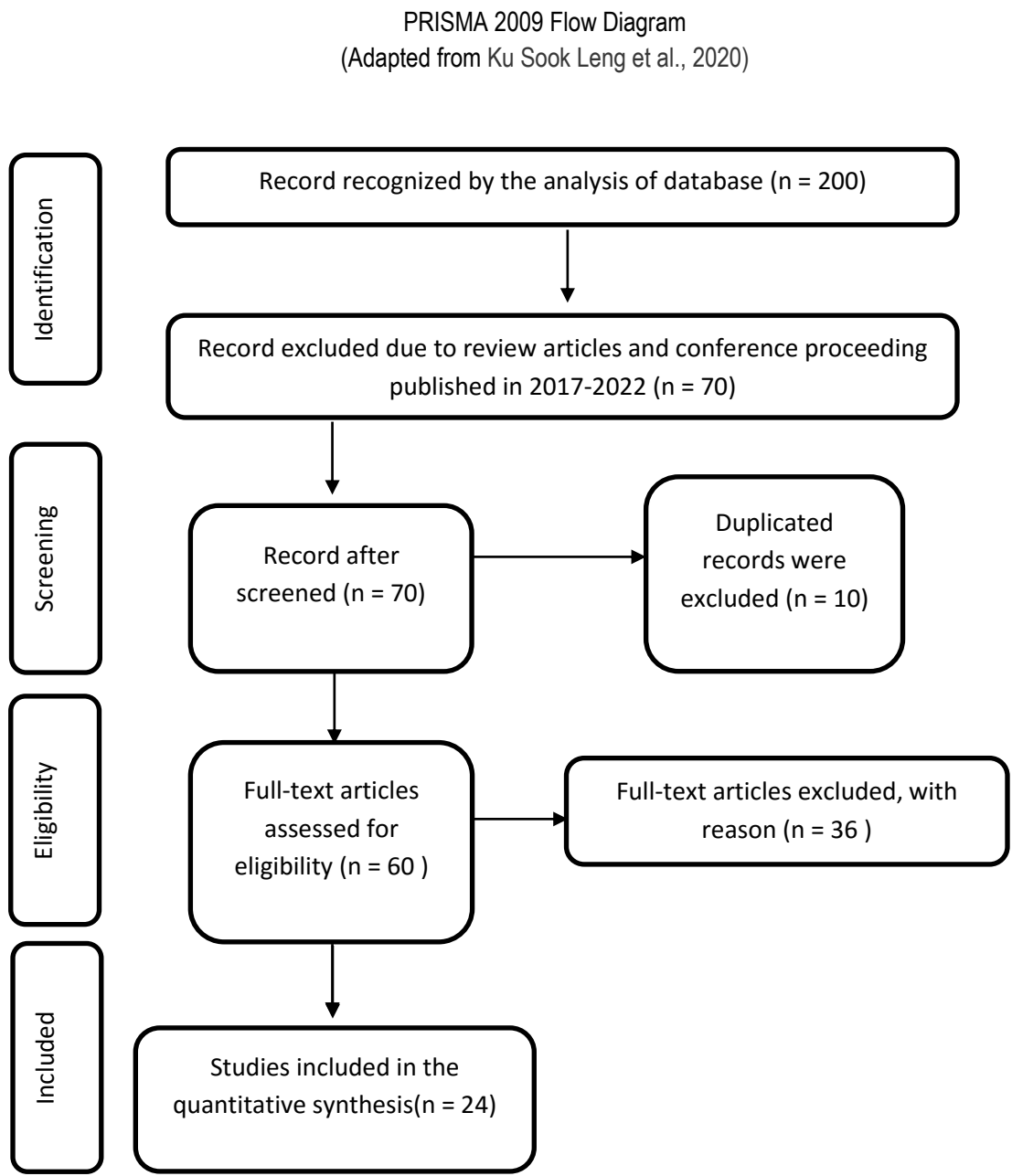

Figure 1. The flow diagram of the study

       The articles that remaining were reviewed and analyzed. The articles papers are reviewed and chosen to be included in the review based on the several criteria (a)ethnomatematics in geometry (b) research in indonesia (c) method (d) conclusion, and (e) symmetry group.

Table 1. Search Results Related to Ethnomathematics in Geometry

| No | Journal/ Proceeding | Method | Conclusion | Symmetry group |
|---|---|---|---|---|
| 1 | Ethnomathematics: Mathematical Values in Masjid Agung Demak (Radiusman, et al, 2021) | Qualitative Research | Masjid Agung Demak contains several mathematical elements, including geometry and transformation. The geometric element is found on the Tajug which consists of three levels, the top of the Tajug, and the location of the Soko Guru pillars. The elements of transformation found in the carving of padma and bledeg. The Masjid Agung Demak is an architectural work that contains elements of mathematics. It does not only contains religious elements but accidentally contains elements of mathematics. | This research only found some elements of mathematics (geometry and transformation) in a small part of Masjid Agung Demak. |
| 2 | Study Ethnomathematics: Revealing Mathematics Ideas on Minangkabau Traditional Weaving Songkets in Pandai Sikek. (Enmufida, Turmudi & Hidayat, A. S., 2020). | Qualitative Research Methods through Ethnographic Approaches. | The result of the research is to reveal a mathematical idea that is used by Pandai Sikek people on manufacturing songket that there are mathematic activities and geometry patterns in traditional songket pattern of Minangkabau. This research enlights people that mathematics has relationship with culture. | There are geometric patterns in the songket fabric pattern is translation, reflection, and rotation |
| 3 | Ethnomathematics: Utilizing South Sumatra's Cultures to Emphasize Prospective Teachers' Creativity in | Descriptive Quantitative Research | Students are able to create more than one mathematical problems, proposes more than one resolution, and some of | The topic engaged with the culture tend to be uniform, with most groups employ the concept of geometry and transformation, namely |
| | Creating Mathematical Problem (Araiku J, Somakim & Pratiwi, W. D., 2019) | | them are able to evaluate those alternatives. | reflection and translation |
| 4 | Ethnomathematics Exploration of the Toba Community: Elements of Geometry Transformation Contained in Gorga (Ornament on Bataks House). (Ditasona, D., 2018) | Qualitative Research | From some gorga motifs can be concluded that the motive has been using the principle of reflection, rotation, translation and dilatation. Gorga made using traditional tools so that in some motives found not too precise. It is interesting to see that the accuracy of their calculations by using simple tools of their day requires a mathematical ability that is good enough for the craftsman in making a motif gorga. | From some gorga motifs can be concluded that the motive has been using the principle of reflection, rotation, translation and dilatation |
| 5 | Learning Geometry and Values From Patterns: Ethnomathematics on The Batik Patterns of Yogyakarta, Indonesia (Prahmana, R.C.I., & D'Ambrosio, U., 2020) | Ethnography Study | Yogyakarta batik uses the concept of geometry transformation in the making of Yogyakarta's unique Batik motif. Besides that, each motif or pattern also contains local values. These, namely moral, historical, and philosophical values, can be felt, reflected, and applied in daily life, such as values that teach leadership, good deeds, and so on. | The people of Yogyakarta have used the concept of geometry transformation in making batik patterns such as in the Babon Angrem, Parang Barong, Parang Klitik, Sidomukti, Semen Bondhat, Sidoluhur, Soblog, and Sidowirasat motifs namely translation, dilatation, reflection and rotation |
| 6 | Buginese Ethnomathematics: Barongko Cake Explorations as Mathematics Learning Resources (Pathuddin, | Qualitative Descriptive with an Ethnographic Approach | This research found that Barongko making process involves mathematics in the concept of division, congruence and similarity, as well as a triangular | This research found that Barongko making process involves mathematics in the concept of division, congruence and similarity, as well as a triangular prism, |

| | | | prism, and half sphere. This cake has the potential to be used as a source of contextual mathematics learning in schools. | and half sphere. This cake has the potential to be used as a source of contextual mathematics learning in schools. |
|---|---|---|---|---|
| 7 | Ethnomathematics: Mathematics in Batik Solo (Faiziyah, N., 2020) | Ethnography Method | Some Batik Solo motifs contain mathematical elements, especially geometric material, namely the principles of translation and reflection. Besides, that there are vertical and horizontal lines as well as perpendicular and parallel lines included. The mathematical concept, especially the geometry used in the making of Batik, is the use of the tessellation or tiling principle. The tessellation principle is the basis for the development of the Solo batik motifs. | This line motif can be used in learning transformation geometry, namely translation, reflection, and the principle of tessellation |
| 8 | Ethnomathematics of the Jami Mosque Jingah River as a Source Mathematics Learning (Fajriah, N and Suryaningsih, Y., 2021) | Ethnographic Approach | Ethnomathematical exploration of the Jami Mosque Jingah River is a mathematical concept of lines, vertices, angles, properties, perimeter and area of a flat shape: rectangle, square, trapezoid, rhombus, circle, half-spherical volume, blank, symmetrical, ellipse, parallel, parallel lines, congruence, condition of the hexagon. | In identifying ecocultural mathematics about space and geometry, four principles were defined and discussed about the structure of language, lines and reference points, size of space and worldview and interpretation of space as a place. |
| 9 | Ethnomathematics: Disclosing Mathematical Concept in Batak Toba | Ethnographic Approach | The mathematical theory consists of a dimensional geometry concept, two-dimensional geometry, | The mathematical theory consists of a dimensional geometry concept, two-dimensional geometry, |
| | Traditional House (Hidayat, E., et al, 2020) | | three-dimensional geometry, geometric transformation, and number patterns. Each of the forms contained in the traditional Toba Batak house contains philosophical values that are used as learning resources. | three-dimensional geometry, geometric transformation, and number patterns. |
| 10 | The students' Mathematics Understanding through Ethnomathematics based on Kejei Dance (Ma'Rifah, N., et al, 2018) | Exploration Studies | The results showed that based on kejei dance, students could build geometry. They have geometry skills through the following activities: identifying and making examples and examples of rejection of flat buildings found in kejei dance. Subjects are able to use models and symbols to represent the concept of building flat on kejei dance. | The mathematical concept applied to kejei dance is to build geometry. Students have geometry skills through the following activities: identifying and making examples and examples of flat-building denials found in kejei dance. |
| 11 | Ethnomathematics: Exploring the Activities of Culture Festival (Maryati and Prahmana, R. C. I., 2018). | Ethnographic Methods | The result shows that in the activities of the culture festival have the concepts of mathematics, such as transformation geometry, cone wall area, cone nets, Pythagorean theorem, and circumference of the circle. The result would be the context in design the learning trajectory by using etnomathematics as the starting point on learning process in realistic mathematics education approach for the future research. | The result shows that in the activities of the culture festival have the concepts of mathematics, such as transformation geometry, cone wall area, cone nets, Pythagorean theorem, and circumference of the circle. |
| 12 | Ethnomathematical | Ethnography | Based on the results and | Toraja's typical carving |
| | Review of Toraja's Typical Carving Design in Geometry Transformation Learning (Nugraha, Y. S., 2018) | Approach | discussion above, it can be concluded that in addition to having philosophical meaning, Toraja's typical carving design also contains mathematical elements, especially in the basic concept of geometry transformation. | design also contains mathematical elements, especially in the basic concept of geometry transformation namely reflection, translation, dilatation, and rotation. |
| 13 | Identifying Geometrical Objects in Sumur Gumuling Tamansari: An Ethnomathematics Analysis (Alvian, D. N., et al, 2020) | Qualitative Method with an Ethnographic Approach | The results of this study are there are mathematical elements (ethnomathematics) in the Sumur Gumuling Tamansari building, namely the concept of a rectangle, pentagon, rhombus, triangle, circle, and the combined concept of a rectangular flat structure with a half-circle. In this case, it can be said that mathematics enters all aspects of life including cultural buildings. | Sumur Gumuling Tamansari building contains mathematical elements in the concept of field geometry, namely the concept of rectangles, pentagons, rhombuses, triangles, circles and the combined concept of rectangular and semicircular flat shapes. |
| 14 | Ethnomathematics on Woven Fabric (Tembe Nggoli) of Mbojo tribe society (Nurbaeti et al, 2018). | Ethnography Method | The results of this study indicate that there is a concept of geometry at woven fabric (tembe nggoli) of mbojo tribe society that can be used as a source of learning and make students and the community better understand the relationship between their culture with the concept of mathematics. | Motif woven fabric of tembe nggoli of mbojo tribes society there are concepts of geometry such as the concept of translation, reflection, and the concept of dilatation as well as the concept of rotation in pattern making nggusu tolu (triangles), nggusu upa pattern (rectangle), waji pado pattern (parallelogram), and the pattern of nggusu waru (octagon). |
| 15 | Ethnomathematics: Geometric Analysis of Historical Buildings Ngawen Temple in | Descriptive Research with a Qualitative | The results of this study were obtained namely, Ngawen Temple set in Buddhism. The Ngawen | The shape of the Ngawen temple resembles the geometry of the cuboid, the rectangular pyramidal |
| | Magelang (Pamungkas, M. D., Zaenuri and Wardono, 2020). | Approach | temple complex consists of five temples that line the parallel from north to south. Temple building facing the east. From the south of Ngawen Temple I, II, III, IV, and V, with each temple to plan the square. One of the uniqueness of Ngawen temple is the existence of 4 lion statues in every corner of temple II and temple IV. | frustum, and the rectangular pyramid. |
| 16 | Exploration of Javanese Culture Ethnomathematics Based on Geometry Perspective (Pramudita, K and Rosnawati, R., 2019) | Exploration | Exploration results show that traditional houses, especially Joglo houses, have ethnomatematics related to lines, angles, triangles, rectangles, quadrilaterals, congruency, and pythagorean theorems. Batik motifs (Lereng, ceplok, and Jlamprang motifs) are ethnomatematics related to lines and angles, triangles, quadrilaterals, circles, and congruency. While the lines, rectangles, circles, and curved side spaces can be associated with ethnomatematics in andong (horse cart). | Traditional houses especially Joglo houses related to lines, angles, triangles, rectangles, quadrilaterals, congruency, and pythagorean theorems. Batik motifs (Lereng, ceplok, and Jlamprang motifs) related to lines and angles, triangles, quadrilaterals, circles, and congruency rectangles, circles, and curved side spaces can be associated with ethnomatematics in andong (horse cart). |
| 17 | Ethnomathematics: Exploration of a Mosque Building and its Ornaments (Purniati, T, Turmudi and Suhaedi, D. (2019) | Qualitative Approach to Exploration Methods | The results of the study show that there are ethnomathematics aspects of the mosque building and its ornaments. The mathematical concepts of the mosque building and its ornaments are related to the concepts of geometry and algebra. Mathematical | The mathematical concepts that can be revealed from the building and ornaments at the Great Mosque of Cimahi are geometry concepts (plane geometry, space geometry, and transformation geometry) and algebra (frieze group). |

| | | | presentations related to mosque buildings and its ornaments are expected to help students to connect mathematics with the surrounding culture so that mathematics becomes more meaningful. | |
|---|---|---|---|---|
| 18 | Ethnomathematics in Balinese Traditional Dance: A Study of Angles in Hand Gestures (Radiusman, R., 2020) | Case study | it can be concluded that the hand gestures of the Pendet Dance can be categorized in three type of angles which are acute angle ($0^0 < \alpha < 90^0$), right angle ($\alpha = 90^0$) and obtuse angle ($90^0 < \alpha < 180^0$). The acute angle can be seen in ngumbang movement, while the right angle can be seen in agem and ulap-ulap movements and the obtuse angle can be seen in the ngelung and ulap-ulap movements. | The hand gestures of the Pendet Dance can be categorized in three type of angles which are acute angle ($0^0 < \alpha < 90^0$), right angle ($\alpha = 90^0$) and obtuse angle ($90^0 < \alpha < 180^0$). |
| 19 | Ethnomathematics: Exploration in Javanese Culture (Risdiyanti, I and Prahmana, R. C. I, 2017) | Ethnography method | The result is exploration ethnomathematics in the several motifs of Yogyakarta batik that contains philosophy, deep cultural value, and mathematics concept, especially geometry transform subject. | Batik of Yogyakarta does applied the geometry transformation concepts such as reflection, translation, rotation, and dilatation. |
| 20 | The Ethnomathematics: Exploration of Gayo Tribe Local Wisdom Related to Mathematics Education (Yustinaningrum, B, Nurliana and | Research and Development with the Four-D Approach. | The traditional measuring tools used in activities of the Gayo community can be developed for teaching measurement topic in the second and third grade. | The traditional measuring tools used in activities of the Gayo can be developed for teaching geometry topics, such as angle, line, plane geometry, and solid |
| | Rahmadhani, E, 2018) | | | geometry. |
| 21 | The Ethnomathematics in Making Woven Bamboo Handicrafts of Osing Community in Banyuwangi, Gintangan Village as Geometry Teaching Material (Yudianto, E et al, 2020). | Ethnography Approach | During the process of calculating the drying time of bamboo and the amount of bamboo whittles, comparative concept of value appeared. During the design process, the concept of quadrilateral and number sequences emerged. The result of this research was geometry teaching material in the form of test package. | The making of these handicrafts patterns were related to the concept of flat shapes and geometric shapes as well as quadrilateral shape. |
| 22 | Integrating Ethnomathematics into Augmented Reality Technology: Exploration, Design, and Implementation in Geometry Learning (Sudirman, S et al, 2019) | The method used ADDIE (Analysis, Design, Development, Implementation, and Evaluation) approach | Ethnomathematics values in the form of uma lengge traditional house, ethnic karo houses, wa rebo traditional house, demak grand mosque, kaghati kolape kite can be used as marker objects on augmented reality technology and to explain geometry concept; The design stages consist of designing integration concepts, AR interfaces, and learning system development. | The applied concept of triangle is a right-angled triangle. The exterior wall is divided into two parts, with the base formation of a square, rectangle, and trapezoid with an upper wall slope of $120^o$. The shape of the drum resembles the concept of a tube in geometry |
| 23 | Ethnomathematics: Mathematical concepts in Tapis Lampung (Susiana, Caswita and Noer, S. H, 2019). | Ethnographic Research | These various motifs turned out to have profound philosophical meanings related to the lives of Lampung people. Geometri concept include One-Dimentional Geometry and Two-Dimentional Geometry, Geometry Transformation includes: concepts of translation, reflection, | Mathematics concept include One-Dimentional Geometry and Two-Dimentional Geometry, Geometry Transformation includes: concepts of reflection, dilation and rotation. |
| | | | dilation and rotation in various types of tapis Lampung including: Tapis cucuk pinggir, tapis jung sarat, tapis limar sekebar, tapis bintang perak. | |
| 24 | Ethnomathematics in Balinese culture as a Learning Material for Logic and Reasoning Geometry (Wijayanti, N.P. A et al, 2018) | Descriptive Qualitative | The results of research show that there are significant mathematical abilities of students to logic and reasoning in solving geometry-based culture of Bali. | The relationship related to the concept of plane figure in mathematics. Mathematics expresses a wide overview of mathematics that includes arithmetic, classification, marketing, mathematics modeling is a culture product. |

From the 24 selected articles, 13 articles used the ethnographic method approach, 5 articles used the qualitative method, 3 articles used the exploratory method, 2 articles used the development method, and 1 article used the case study method. Ethnographic research is currently becoming a study with many enthusiasts, especially in researches in the social sciences and humanities as well as in research related to ethnicity and culture. This ethnographic method is widely used to examine people's lives obtained from community members based on the adaptation patterns possessed by researchers (Hanifah, 2010). Ethnography as a method is used as a tool to explain in detail about the culture or behavior patterns of a society.

The geometric concepts that exist in batik and weaving or fabric motifs, Yogyakarta batik uses the concept of geometry transformation in the making of Yogyakarta's unique Batik motif (Babon Angrem, Parang Baron, Parang Klitik, Sidomukti, Semen Bondhat, Sidoluhur, Soblog, and Sidowirasat) and Bima (Tembe Nggoli) weaving, namely translation, dilation, reflection and rotation. In addition to the four concepts above, in Solo batik there is an addition to the concept of geometry, namely tessellation. In Tapis Lampung there are only 3 geometric concepts, namely reflection, dilation and rotation. Meanwhile, in Minangkabau traditional weaving Songkets in Pandai Sikek there are geometric concepts, namely translation, reflection and rotation.

In historic buildings such as the Great Mosque of Demak, the geometric element is found on the Tajug which consists of three levels, the top of the Tajug, and the location of the Soko Guru pillars. The elements of transformation found in the carving of padma and bledeg. In addition, other historical buildings, namely the Jami Mosque Jingah River, Batak Toba Traditional House, Batak Traditional House (Ruma Gorga) and carving designs of Toraja traditional houses found the concepts of translation, dilation, reflection and rotation. The mathematical concepts that can be revealed from the building and ornaments at the Great Mosque of Cimahi are geometry concepts (plane geometry, space geometry, and transformation geometry) and algebra (frieze group).

At Ngawen Temple Magelang, the shape of the Ngawen temple resembles the geometry of the cuboid, the rectangular pyramidal frustum, and the rectangular pyramid. Traditional houses especially Joglo houses related to lines, angles, triangles, rectangles, quadrilaterals, congruency, and pythagorean theorems. Meanwhile, Sumur Gumuling Tamansari building contains mathematical elements in the concept of field geometry, namely the concept of rectangles, pentagons, rhombuses, triangles, circles and the combined concept of rectangular and semicircular flat shapes.

By integrating Augmanted Reality Technology, ethnomathematics values in the form of uma lengge traditional house, ethnic karo houses, wa rebo traditional house, demak grand mosque, kaghati kolape kite can be used as marker objects on augmented reality technology and to explain geometry concept. The applied concept of triangle is a right-angled triangle. The exterior wall is divided into two parts, with the base formation of a square, rectangle, and trapezoid with an upper wall slope of $120^0$. The shape of the drum resembles the concept of a tube in geometry.

In the traditional dances from Bengkulu, namely the Kejei Dance and the Balinese Traditional Dance, namely the Pendet Dance, the concept of geometry is also found. The mathematical concept applied to kejei dance is to build geometry. Students have geometry skills through the following activities: identifying and making examples and examples of flat-building denials found in kejei dance. While in Balinese Traditional Dance, the hand gestures of the Pendet Dance can be categorized in three type of angles which are acute angle $(0^0 < \alpha < 90^0)$, right angle $(\alpha = 90^0)$ and obtuse angle. The acute angle can be seen in ngumbang movement, while the right angle can be seen in agem and ulap-ulap movements and the obtuse angle can be seen in the ngelung and ulap-ulap movements. The angle construction is performed by the gestures of the elbow, upper and lower arm, and their position towards other parts of the body such as the chest and eyes.

## CONCLUSION

From the findings and the results of the articles included, it is clear that the results of this study indicate that ethnomathematics in many cultures in Indonesia uses the concept of geometry, including mosques, temples, batik or woven fabrics, typical food of an historical place, even in regional dances. It was recommended that educators make use of cultural examples that are inclusive of all students coming from a variety of cultural backgrounds. Learners could also be encouraged by their teachers to look for connections between the ideas of geometry they are studying and instances from their own culture and environment.